\documentclass{amsart}

\title{Evidence for The Gromov-Witten/Donaldson-Thomas Correspondence}
\author{Amin Gholampour}

\address{
Department of Mathematics\\
University of British Columbia \\
1984 Mathematics Road \\
Vancouver, B.C., Canada V6T 1Z2}
 
\email{amin@math.ubc.ca}

\author{Yinan Song}
\address{
Department of Mathematics\\
University of British Columbia \\
1984 Mathematics Road \\
Vancouver, B.C., Canada V6T 1Z2}
\email{yinan@math.ubc.ca}

\date{September 30, 2005}
\usepackage{amsmath,amsthm,amsfonts,amscd}

\newtheorem{thm}{Theorem}
\newtheorem{lem}[thm]{Lemma}
\newtheorem{rem}[thm]{Remark}

\newtheorem{prop}[thm]{Proposition}

\newcommand{\Norm}{\operatorname{Norm}}

\newcommand{\Def}{\operatorname{Def}}
\newcommand{\Ob}{\operatorname{Ob}}
\newcommand{\id}{\operatorname{id}}
\newcommand{\Hilb}{\operatorname{Hilb}}

\begin{document}

\maketitle
\begin{abstract}
We study the equivariant Gromov-Witten and Donaldson-Thomas
theories of $\mathbf{P}^2$-bundles over curves. We show the  equivariant GW/DT correspondence holds to first
order for certain curve classes.
\end{abstract}
\section{Introduction}  \label{sec:intro}

Both Gromov-Witten and Donaldson-Thomas theory arise from the
enumerative geometry of threefolds. In \cite{MNOP1}, \cite{MNOP2},
and \cite{Bryan-Pandharipande}, correspondences between these two
theories were conjectured in various settings. In this paper, we
provide evidence for the equivariant conjecture in the case
where the threefold is a $\mathbf{P}^2$-bundle over a curve. The
equivariant Gromov-Witten theory of $\mathbf{P}^2$-bundle over
curves is solved in \cite{Amin} for certain curve classes, but we
can only compute the Donaldson-Thomas invariants up to first
order. Still, our computation is the first compact non-toric example where
the correspondence has been checked.

Let $C$ be a smooth and connected curve of genus $g$. Let $L_1$ and
$L_2$ be line bundles of degrees $k_1$ and $k_2$, and let $X$ be the
total space of the $\mathbf{P}^2$-bundle $\mathbf{P}(\mathcal{O}
\oplus L_1 \oplus L_2)$ over $C$. The three dimensional torus,
$(\mathbf{C}^*)^3$ acts diagonally on the fibers of $X$. We
denote by $t_0$, $t_1$ and $t_2$, the equivariant parameters of
this action along the first, second and third summand,
respectively. Let $\beta$ be a curve class in $X$.

Following \cite{Amin}, we define the partition function for the
degree $\beta$ equivariant Gromov-Witten invariants of $X$ as:

\begin{align*}
Z^{GW}_{\beta}(g\operatorname{|}k_{1},k_{2})&=\sum_{h=0}^{\infty}u^{2h-2-K_{X}\cdot\,
\beta}\int_{[\overline{M}_{h}^{\bullet}(X,\beta)]^{vir}}1,
\end{align*} where $\overline{M}_{h}^{\bullet}(X,\beta)$ is the moduli
space of the genus $h$, degree $\beta$, possibly disconnected
stable maps to $X$, and $K_{X}$ is the canonical class of $X$. The
integration in the above formula denotes equivariant pushforward
to a point, so it can be non-zero for negative virtual dimension.
The Gromov-Witten partition function will be an element in
\[\mathbf{Q}((u))[t_0, t_1, t_2].\]

As in \cite{Bryan-Pandharipande}, we define the the partition function for the
class $\beta$ equivariant Donaldson-Thomas invariants of $X$ as

\[
Z_{\beta}^{DT}(g\operatorname{|}k_{1},k_{2})=\sum_{n\in\mathbb{Z}}q^{n}\int_{[I_{n}(X,\beta)]^{vir}}1,
\]
where $I_n(X, \beta)$ is the moduli space of ideal sheaves
parameterizing subschemes $Y$ in $X$, whose maximal one
dimensional component is in the class $\beta$ and $\chi(O_{Y})=n$.
Again, the integral denotes equivariant pushforward to a point.
The Donaldson-Thomas partition function is an element in
\[\mathbf{Z}((u))[t_0, t_1, t_2].\]
It takes coefficients in $\mathbf{Z}((u))$, rather than
$\mathbf{Q}((u))$ because $I_n(X, \beta)$ is a scheme, so the
virtual class is defined over $\mathbf{Z}$.

We define the reduced Donaldson-Thomas partition function by
$$Z_{\beta}^{DT}(g\operatorname{|}k_{1},k_{2})'=M(-q)^{-\int_{X} c_{3}(T_{X}\otimes
K_{X})}\cdot Z_{\beta}^{DT}(g\operatorname{|}k_{1},k_{2}),$$ where
$$M(q)=\prod_{n \geq 1}\frac{1}{(1-q^n)^n}$$ is the McMahon
function, and $T_{X}$ is the tangent bundle of $X$ (See Definition
9.2, and Conjecture 1, in\cite{Bryan-Pandharipande}).

In this notation, GW/DT correspondence conjecture in
\cite{Bryan-Pandharipande} reads as follows:
\subsection*{GW/DT Conjecture.} After the change of variable,
$e^{iu}=-q$, we have
\begin{equation} \label{equ:GW/DT}
(-i)^{-K_{X}\cdot\,\beta}Z^{GW}_{\beta}(g\operatorname{|}k_{1},k_{2})=(-q)^{\frac{1}{2}K_{X}\cdot
\,\beta}Z_{\beta}^{DT}(g\operatorname{|}k_{1},k_{2})'.
\end{equation}

We define $\beta_0$ to be the curve class of the section given by the locus $(1: 0:
0)$, and $f$ to be the line class of a fibre in $X$. Then we prove
the following result:

\begin{thm} \label{result}
GW/DT Conjecture holds for the leading term in $q$, when $\beta=\beta_{0}$
and $\beta=\beta_{0}+f$, and $k_{1},k_{2}\leq0$.
\end{thm}

The case $\beta=\beta_{0}$ is straightforward. We will prove the
theorem for $\beta=\beta_{0}+f$, and when $k_{1},k_{2} < -1$,
where we can express our results more uniformly. The case, where
one or both of $k_{1}$ and $k_{2}$ are equal to -1 or 0, is
similar. In Section 2, we will compute the left hand side of
(\ref{equ:GW/DT}) (Proposition \ref{GW}). In section 3, we will
compute the first term of the series in the right hand side of
(\ref{equ:GW/DT}) (Proposition \ref{DT}). Theorem~\ref{result}
will follow after changing the variable, $e^{iu}=-q$, and
comparing the coefficients of $q$ with the lowest power in the
both sides of (\ref{equ:GW/DT}).

\section{The Gromov-Witten Theory of $X$}
Let $X$ and $\beta=\beta_{0}+f$ be as above. We can express the result of
this section as follows:

\begin{prop} \label{GW} Let $k_{1}, k_{2}< -1$, then the partition function for the degree $\beta=\beta_{0}+f$
equivariant Gromov-Witten invariants of $X$ is given by
\begin{align*}&Z_{\beta}(g\operatorname{|}k_{1},k_{2})=
(t_{0}-t_{1})^{g-k_{1}-3}(t_{0}-t_{2})^{g-k_{2}-3}\\
&\hspace{.5cm}\cdot\big((4g-4-k_{1}-k_{2})t_{0}-(2g-2-k_{2})t_{1}-(2g-2-k_{1})t_{2})\big)
\big(2\operatorname{sin}\frac{u}{2}\big)^{k_{1}+k_{2}+3}.
\end{align*}
\end{prop}

In \cite{Amin}, the section class equivariant Gromov-Witten
theory of $X$ was completely determined. The following matrices
were defined in \cite{Amin}:

$$ A=
\left[\begin{array}{ccc} (t_{0}-t_{1})(t_{0}-t_{2}) & 0
& 0 \\
0 & (t_{1}-t_{0})(t_{1}-t_{2})
& 0 \\
0 & 0 & (t_{2}-t_{0})(t_{2}-t_{1})
\end{array}
\right], $$
$$
B= \left[\begin{array}{ccc}
\frac{2(2t_{0}-t_{1}-t_{2})}{(t_{0}-t_{1})(t_{0}-t_{2})} &
\frac{t_{0}+t_{1}-2t_{2}}{(t_{0}-t_{1})(t_{0}-t_{2})}
&\frac{t_{0}+t_{2}-2t_{1}}{(t_{0}-t_{1})(t_{0}-t_{2})} \\
\frac{t_{0}+t_{1}-2t_{2}}{(t_{1}-t_{0})(t_{1}-t_{2})} &
\frac{2(2t_{1}-t_{0}-t_{2})}{(t_{1}-t_{0})(t_{1}-t_{2})}
&\frac{t_{1}+t_{2}-2t_{0}}{(t_{1}-t_{0})(t_{1}-t_{2})} \\
\frac{t_{0}+t_{2}-2t_{1}}{(t_{2}-t_{0})(t_{2}-t_{1})} &
\frac{t_{1}+t_{2}-2t_{0}}{(t_{2}-t_{0})(t_{2}-t_{1})} &
\frac{2(2t_{2}+t_{0}-t_{1})}{(t_{2}-t_{0})(t_{2}-t_{1})}
\end{array} \right]\phi^{3}.
$$

$$
M_{1}=\left[\begin{array}{ccc}
t_{0}-t_{1} & 0 & 0 \\
0 &
0 & 0 \\
0 & 0 & t_{2}-t_{1}
\end{array} \right]\phi^{-1}, \quad
M_{2}=\left[
\begin{array}{ccc}
t_{0}-t_{2} & 0 & 0 \\
0 &
t_{1}-t_{2} & 0 \\
0 & 0 & 0
\end{array}\right]\phi^{-1},$$
$$
N=\left[\begin{array}{ccc} \frac{1}{(t_{0}-t_{1})(t_{0}-t_{2})} &
\frac{1}{(t_{0}-t_{1})(t_{0}-t_{2})}
& \frac{1}{(t_{0}-t_{1})(t_{0}-t_{2})} \\
\frac{1}{(t_{1}-t_{0})(t_{1}-t_{2})} &
\frac{1}{(t_{1}-t_{0})(t_{1}-t_{2})}
& \frac{1}{(t_{1}-t_{0})(t_{1}-t_{2})} \\
\frac{1}{(t_{2}-t_{0})(t_{2}-t_{1})} &
\frac{1}{(t_{2}-t_{0})(t_{2}-t_{1})} &
\frac{1}{(t_{2}-t_{0})(t_{2}-t_{1})}
\end{array} \right]\phi^{2},
$$
where $\phi=2\operatorname{sin}\frac{u}{2}$.

In \cite{Amin}, $G=A+B$ was called the genus adding operator, and
$U^{-1}_{1}=M_{1}+N$ and $U_{2}^{-1}=M_{2}+N$ were called the
first and the second level annihilation operator, respectively.

\subsection*{Notation.} For two matrices $U$ and $V$, by
$(U^{a},V^{b})$ for $a,b \in \mathbb{Z}^{+}$, we mean the sum of
the all products that we can write containing $a$ copies of $U$
and $b$ copies of $V$. For example
$$(U^{2},V)=U^{2}V+UVU+VU^{2}.$$

Now by Theorem 1.7 in \cite{Amin}, for $k_{1}, k_{2}\leq 0$, and
$r \geq 0$, we have
$$Z_{\beta_{0}+rf}(g\operatorname{|}k_{1},k_{2})=\sum\operatorname{tr}
\left((A^{a},B^{b})(M_{1}^{m_{1}},N^{n_{1}})(M_{2}^{m_{2}},N^{n_{2}})\right),$$
where the sum is over all nonnegative integers
$a$,\hspace{1mm}$b$,\hspace{1mm}$m_{1}$,\hspace{1mm}$m_{2}$,\hspace{1mm}$n_{1}$,
and $n_{2}$ such that
$$\begin{array}{llll} a+b=g-1, & m_{1}+n_{1}=-k_{1}, &
m_{2}+n_{2}=-k_{2}, & b+n_{1}+n_{2}=r
\end{array}.$$

\textsc{Proof of Proposition \ref{GW}:} By the discussion above,
for $\beta=\beta_0+f$, we have
\begin{align*}
&Z_{\beta}(g\operatorname{|}k_{1},k_{2})=\\&\operatorname{tr}\left((A^{g-2},B)M_{1}^{-k_{1}}M_{2}^{-k_{2}}+
A^{g-1}(M_{1}^{-k_{1}-1},N)M_{2}^{-k_{2}}+A^{g-1}M_{1}^{-k_{1}}(M_{2}^{-k_{2}-1},N)\right).
\end{align*}
Now since $\operatorname{tr}(UV)=\operatorname{tr}(VU)$ for any
two matrices, and also since $A$, $M_{1}$, and $M_{2}$ mutually
commute, this is equal to
\begin{align*}
&(g-1)\operatorname{tr}\left(A^{g-2}M_{1}^{-k_{1}}M_{2}^{-k_{2}}B\right)
-k_{1}\operatorname{tr}\left(A^{g-1}M_{1}^{-k_{1}-1}M_{2}^{-k_{2}}N\right)\\
&-k_{2}\operatorname{tr}\left(A^{g-1}M_{1}^{-k_{1}}M_{2}^{-k_{2}-1}N\right)
.\end{align*}

From this, Proposition~\ref{GW} easily follows. \qed

\section{the Donaldson-Thomas Theory of X}

Let $X$ and $I_n(X, \beta)$ be as in Section \ref{sec:intro}, and let
$\pi$ be the projection from $X$ to $C$. Unlike the computation in
Gromov-Witten theory, we can only compute the first term of the
partition function in Donaldson-Thomas theory in the right hand
side of (\ref{equ:GW/DT}). The power of $q$ of the first term is
the holomorphic Euler characteristic of the structure sheaf of a
purely 1-dimensional subscheme, $Y$, in the class $\beta$, which
is equal to $\chi(\mathcal{O}_Y)=1-g$. From now on we fix
\[n=1-g.\]

By definition, the coefficient of the first term of the right hand
side of (\ref{equ:GW/DT}) is \begin{equation} \label{equ:def
N}N^{DT}_{n}(X, \beta)=\int_{[I_n(X, \beta)]^{vir}}1 =
 \int_{[I_n(X,
 \beta)^{T}]^{vir}}\frac{1}{e(\Norm^{vir})},\end{equation}
where $I_n(X, \beta)^{T}$ is the fixed loci of the induced torus
action on $I_n(X, \beta)$, and $\Norm^{vir}$ is its virtual normal
bundle, and $e$ denotes equivariant Euler class, see \cite{Bryan-Pandharipande} for a detailed discussion.
Now we express the result of this section:

\begin{prop} \label{DT}
Let $X$, $n$, and $k_1$, $k_2 < -1$ be as above. Then the
equivariant Donaldson-Thomas invariant of class $\beta = \beta_{0}+f$ and
holomorphic Euler characteristic $n$ is given by:
\begin{align*}
N^{DT}_n(X,\beta)=&
(t_{0}-t_{1})^{g-k_{1}-3}(t_{0}-t_{2})^{g-k_{2}-3}\\
&\cdot\big((4g-4-k_{1}-k_{2})t_{0}-(2g-2-k_{2})t_{1}-(2g-2-k_{1})t_{2})\big).
\end{align*}
\end{prop}
Before giving the proof, we first describe the fixed loci
explicitly:

\begin{lem}
$I_n(X, \beta)^{T}$ is isomorphic to the union of two disjoint
copies of $C$, and moreover the universal subscheme in $I_n(X,
\beta)^{T} \times X$ is a complete intersection.
\end{lem}

\textsc{Proof:}
Let $M = C_1 \coprod C_2$ where $C_1$ and $C_2$ are two copies of $C$. We will consider $M$ as $\Hilb^1(M)$ and show that $M$ is isomorphic to $I_n(X, \beta)^T$ by showing that they represent isomorphic functors of points. In other words, we will show how to construct a family of one-point in M from a family of fixed subschemes and vice versa. First we construct an $M$-family of fixed subschemes:
\[\mathcal{Y} \subset M \times X \stackrel{pr_1}{\rightarrow} M,\]
where $pr_1$ and $pr_2$ denote projections from $M \times X$ onto $M$ and $X$ respectively. We define $\Delta$ to be the diagonal divisor in $C \times C$. We use $\Delta_1$ and $\Delta_2$ to distinguish between the two diagonals in $M \times C$. We use $\id \times \pi$ denote from the map from $M \times X$ to $M \times C$.

The projective bundle $X$ can be constructed in three ways, as $\mathbf{P}(\mathcal{O}\oplus L_1 \oplus L_2)$, and $\mathbf{P}(L_1^{*}\oplus \mathcal{O} \oplus L_2 L_1^{*})$, and $\mathbf{P}(L_2^{*} \oplus L_1 L_2^{*} \oplus \mathcal{O})$. The tautological bundles change depending on the construction. We use the following notation for the hyperplane divisors:

\begin{align*}
H_0= \mathbf{P}(L_1 \oplus L_2) &\subset X = \mathbf{P}(\mathcal{O} \oplus L_1 \oplus L_2) \\
H_1= \mathbf{P}(L_1^{*} \oplus L_2L_1^{*}) &\subset X = \mathbf{P}(L_1^{*} \oplus \mathcal{O} \oplus L_2 L_1^{*}) \\
H_2= \mathbf{P}(L_2^{*} \oplus L_1 L_2^{*}) &\subset X =
\mathbf{P}(L_2^{*} \oplus L_1 L_2^{*} \oplus \mathcal{O})
\end{align*}

The corresponding line bundles are $\mathcal{O}_X(1)$, $\pi^*L_1 \otimes \mathcal{O}_{X}(1)$ and $\pi^*L_2 \otimes \mathcal{O}_{X}(1)$. Finally, we let $D_1 = pr_2^*H_1$ and $D_2=pr_2^*H_2 \cup (\id \times \pi)^*\Delta_1$ and $D'_{1} = pr_2^*H_2$ and $D'_2=pr_2^*H_1 \cup (\id \times \pi)^*\Delta_2$.

The family $\mathcal{Y}$ is a disjoint union of two components $\mathcal{Y}_1$ and $\mathcal{Y}_2$. One component is the intersection of $D_1$ and $D_2$ and the other is the intersection of $D'_1$ and $D'_2$. This family is fixed under the $T$-action.

This family $\mathcal{Y}$ is constructed out of the universal family of $\Hilb^1(M)$. Therefore, given a $B$-family of one-point in $M$, which is equivalent to a map from $B$ to $M$, we can pull back the family $\mathcal{Y}$ to obtain a family of $T$-fixed subschemes over $B$.

Conversely, let's consider a $B$-family $\mathcal{Z}$ of fixed subschemes. We can intersect them with the hyperplane divisor at the infinity, $H_0$. Since the subschemes are of the class $\beta_0 + f$, each fixed scheme intersects $H_0$ at a point. This point is a $T$-fixed point, so it is a point on either the locus of $(0:1:0)$ in $X = \mathbf{P}(L_1^{*} \oplus \mathcal{O} \oplus L_2 L_1^{*})$ or the locus of $(0:0:1)$ in $X = \mathbf{P}(L_2^{*} \oplus L_1 L_2^{*} \oplus \mathcal{O})$. The union of these two loci is isomorphic to $M$. Therefore, we get a $B$-family of one-points in $M$. It is easy to see that it is an isomorphism of functors, so we get an isomorphism of the moduli spaces.

\qed

 \textsc{Proof of Proposition \ref{DT}:}
Consider the first component constructed in the previous
lemma. For convenience, we will suppress the subscript in $C_1$ and $\mathcal{Y}_1$ from now on. Recall that $\pi$ and $\id \times \pi$ are the projections from X and $C \times X$ to C and $C \times C$ respectively, that $pr_1$ and $pr_2$ are projections from $C \times X$ onto $C$ and $X$ respectively, and that $p_1$ and $p_2$ are projections from $C \times C$ onto the first and the second factor respectively. Recall that $\mathcal{Y}$ is a zero section of rank two vector bundle on $C \times X$, which is the direct sum of the following two line bundles:

\begin{align*}
\mathcal{O}_{C \times X}(D_1) &= pr_2^{*}(\pi^{*}L_1 \otimes \mathcal{O}_{X}(1)) \\
\mathcal{O}_{C \times X}(D_2) &= pr_2^{*}(\pi^{*}L_2 \otimes \mathcal{O}_{X}(1)) \otimes (\id \times \pi)^{*}\mathcal{O}_{C \times C}(\Delta).
\end{align*}

The K-theory class of the perfect obstruction theory on $C$, considered as
a connected component of $I_n(X, \beta)^T$, is given by
\[\mathbf{R}^{\bullet}pr_{1,*}(\mathcal{H}om(\mathcal{I},
\mathcal{I})-\mathcal{O}_{C \times X}),\] where $\mathcal{I}$ is
the universal ideal sheaf (c.f. ~\cite{MNOP1} page 19).  According to \cite{Gr-Pa}, in terms of
equivariant K-theory classes, we can regard
\[\mathbf{R}^{\bullet}pr_{1,*}(\mathcal{H}om(\mathcal{I}, \mathcal{I})-\mathcal{O}_{C \times X})=
-\Def^{m}-\Def^f+\Ob^m+\Ob^f,\] where the superscripts, $f$, and
$m$, mean the fixed, and the moving parts under the induced action of the
torus, respectively. Since the fixed loci is smooth, we will get
\[N^{DT}_{n}(X, \beta)=\int_{[I_n(X, \beta)^{T}]}\frac{e(\Ob^{f})e(\Ob^{m})}{e(\Def^{m})}.\]
We will show later that $\Ob^{f}$ is zero, so we can write
\begin{equation} \label{eqn:DT}
N^{DT}_{n}(X, \beta)=\int_{C}e(\mathbf{R}^{\bullet}pr_{1,*}(\mathcal{H}om(\mathcal{I},
\mathcal{I})-\mathcal{O}_{C \times X}))^m.
\end{equation}

Since $\mathcal{Y}$ is the complete intersection of two divisors,
$D_1$ and $D_2$, we have the standard Kozsul resolution of the ideal sheaf of $\mathcal{Y}$:
\[0 \rightarrow \mathcal{O}_{C \times X}(-D_1-D_2) \rightarrow \mathcal{O}_{C\times X}(-D_1) \oplus \mathcal{O}_{C \times X}(-D_2) \rightarrow \mathcal{I} \rightarrow 0.\]
Therefore, the equivariant K-theoretic class of $\mathcal{I}$ is
\[\mathcal{O}_{C \times X}(-D_1)+\mathcal{O}_{C \times X}(-D_2) - \mathcal{O}_{C \times X}(-D_1-D_2).\]
Then a formal calculation of equivariant K-theoretic classes gives
us
\begin{align}  \label{align:defobclass}
\mathbf{R}^{\bullet}pr_{1,*}(\mathcal{H}om(\mathcal{I},
\mathcal{I})-\mathcal{O}_{C \times X})
&=-\mathbf{R}^{\bullet}pr_{1,*}(\mathcal{O}_{C \times X}(D_1))-
\mathbf{R}^{\bullet}pr_{1,*}(\mathcal{O}_{C \times X}(D_2))\nonumber\\
&\hspace{.3cm}-\mathbf{R}^{\bullet}pr_{1,*}(\mathcal{O}_{C \times
X}(-D_1))-\mathbf{R}^{\bullet}pr_{1,*}(\mathcal{O}_{C \times X}(-D_2)\nonumber\\
&\hspace{.3cm}+\mathbf{R}^{\bullet}pr_{1,*}(\mathcal{O}_{C \times
X}(D_1-D_2)+2\mathbf{R}^{\bullet}pr_{1,*}(\mathcal{O}_{C \times
X})\\
&\hspace{.3cm}+ \mathbf{R}^{\bullet}pr_{1,*}(\mathcal{O}_{C \times
X}(D_2-D_1)). \nonumber
\end{align}
The map, $pr_1$, can be factored through
\[C \times X \stackrel{\id \times \pi}{\rightarrow} C \times C \stackrel{p_1}{\rightarrow} C,\]
where $p_1$ is the projection onto the first factor, and therefore
\[\mathbf{R}^{\bullet}pr_{1,*} = \mathbf{R}^{\bullet}p_{1,*} \circ \mathbf{R}^{\bullet}(\id \times \pi)_{*}.\]

We will use the following facts in the rest of the proof and we summarize them here in the following lemma:
\begin{lem} \label{facts}
With the same notation as in Proposition \ref{DT}, we have the
following identities.
\begin{align*}
&\mathbf{R}^{\bullet}(\id \times \pi)_{*}(\mathcal{O}_{C \times X}) =
\mathcal{O}_{C \times C}, \nonumber\\
&\mathbf{R}^{\bullet}(\id \times \pi)_{*}(pr_2^{*}\mathcal{O}_{X}(1)) =
p_2^*L_0^{*} \oplus p_2^*L_1^{*}\oplus p_2^*L_2^{*},\nonumber\\
&\mathbf{R}^{\bullet}(\id \times \pi)_{*}(pr_2^{*}\mathcal{O}_{X}(-1)) = 0,\nonumber\\
&\mathbf{R}^{\bullet} p_{1,*}(p_2^* (L_0^i L_1^j L_2^l)) = (1-g+j k_1 +
l k_2) \mathbf{C}_{it_0+jt_1+lt_3},\\
&\mathbf{R}^{\bullet} p_{1,*}(p_2^* (L_0^i L_1^j L_2^l) \otimes \mathcal{O}_{C \times C} (-\Delta) = (1-g+j
k_1+l k_2) \mathbf{C}_{i t_0 + j t_1+ l t_2} - L_0^i L_1^j
L_2^l,\nonumber\\
&\mathbf{R}^{\bullet} p_{1,*}(p_2^* (L_0^i L_1^j L_2^l) \otimes \mathcal{O}_{C \times C}(\Delta)) = (1-g+j
k_1 + l k_2) \mathbf{C}_{i t_0 + j t_1 + l t_2} + L_0^i L_1^j
L_2^l T_C,\nonumber
\end{align*}
where $\mathbf{C}_{s}$ denotes the trivial line bundle on $C$ with
weight $-s$, and $L_0$ denotes the trivial line bundle on $C
\times C$ with weight $-t_0$.
\end{lem}

\textsc{Proof:} The first three are elementary. The fourth follows
from Riemann-Roch. The last two are proved by applying
$\mathbf{R}^{\bullet}p_{1, *}$ to the divisor sequence of $\Delta$
twisted by $p_2^*(L_0^{i}L_1^{j}L_2^{l})$ and
$p_2^*(L_0^{i}L_1^{j}L_2^{l})\otimes \mathcal{O}_{C \times
C}(\Delta)$ respectively. \qed

Now by the first and fourth equalities in Lemma \ref{facts}, we
have
\[\mathbf{R}^{\bullet}pr_{1,*}(\mathcal{O}_{C \times X}) = (1-g) \mathbf{C}_{0}.\]
The second and fourth equalities in Lemma \ref{facts}, implies
that
\begin{align*}
\mathbf{R}^{\bullet}pr_{1,*}(pr_2^*(\mathcal{O}_{X}(1) \otimes \pi^* L_1))
&=\mathbf{R}^{\bullet}(p_1)_{*}((p_2^*L_0^{*} \oplus p_2^*L_1^{*} \oplus p_2^*L_2^{*}) \otimes p_2^*L_1) \\
&=\mathbf{R}^{\bullet}(p_1)_{*}(p_2^*(L_0^{*}L_1) \oplus \mathcal{O} \oplus p_2^*(L_2^{*} L_1)) \\
&=(k_1+1-g) \mathbf{C}_{t_1-t_0}+(1-g) \mathbf{C}_{0} \\
&\hspace{3mm}+(k_1-k_2 +1-g) \mathbf{C}_{t_1 - t_2},
\end{align*}
and similarly, using the second and the sixth equalities in
 Lemma \ref{facts}, we have
\begin{align*}
\mathbf{R}^{\bullet} pr_{1,*}(pr_2^*(\mathcal{O}_{X}(1)& \otimes
\pi^* L_2) \otimes (\id \times \pi)^* (\mathcal{O}_{C \times C}
(\Delta)))=\\
&=\mathbf{R}^{\bullet} p_{1,*}((p_2^*L_0^{*} \oplus p_2^*L_1^{*}
\oplus p_2^*L_2^{*})
\otimes p_2^*L_2 \otimes \mathcal{O}_{C \times C}(\Delta))\\
&=\mathbf{R}^{\bullet} p_{1,*}(p_2^*(L_0^{*} L_2) (\Delta) \oplus
p_2^*(L_1^*
L_2) (\Delta) \oplus \mathcal{O}_{C \times C}(\Delta))\\
&=(1-g+k_2)\mathbf{C}_{t_2-t_0}+ L_0^* L_2 T_C +
(1-g) \mathbf{C}_{0}\\
&\hspace{3mm}+(1-g-k_1+k_2)\mathbf{C}_{-t_1+t_2}+L_1^* L_2 T_C +
T_C.
\end{align*}
Finally, we can write
\begin{align*}
\mathbf{R}^{\bullet}pr_{1,*}((\id \times \pi)^*(p_2^*(L_1^* L_2) (\Delta)))
&=\mathbf{R}^{\bullet} (p_1)_*(p_2^*(L_1^* L_2) (\Delta)) \\
&=(1-g -k_1+k_2) \mathbf{C}_{-t_1+t_2} + L_1^* L_2 T_C,
\end{align*}
and
\begin{align*}
\mathbf{R}^{\bullet} pr_{1,*} ((\id \times \pi)^*(p_2^*(L_1 L_2^* )(-\Delta)))
&=\mathbf{R}^{\bullet} (p_1)_{*} (p_2^*(L_1 L_2^*) (-\Delta)) \\
&=(1-g +k_1-k_2) \mathbf{C}_{t_1-t_2} - L_1 L_2^*,
\end{align*}
where for the first one, we used first and the sixth equalities in
Lemma \ref{facts}, and for the second one, we used the first and
fifth equalities in Lemma \ref{facts}.

So far, we have computed all the terms in the right hand side of
(\ref{align:defobclass}), so we can write
\begin{align} \label{equ:1}
e(\mathbf{R}^{\bullet}pr_{1,*}(\mathcal{H}om(\mathcal{I},
\mathcal{I})- \mathcal{O}_{C \times X}&))^m = \nonumber\\
&\frac{(t_0-t_1)^{-k_1+g-1}(t_0-t_2)^{-k_2+g-1}}{((t_0-t_2)+(k_2+2-2g)[p\,])((t_2-t_1)+(k_1-k_2)[p\,])},
\end{align}
where $[p\,]$ is the class of a point in $C$.

Similar computations for the other component of $I_n(X, \beta)^T$,
yield
\begin{align} \label{equ:2}
e(\mathbf{R}^{\bullet}pr_{1,*}(\mathcal{H}om(\mathcal{I},
\mathcal{I})- \mathcal{O}_{C \times X}&))^m =\nonumber\\
&\frac{(t_0-t_1)^{-k_1+g-1}(t_0-t_2)^{-k_2+g-1}}{((t_0-t_1)+(k_1+2-2g)[p\,])((t_1-t_2)+(k_2-k_1)[p\,])}.
\end{align}

Equation (\ref{align:defobclass}) also shows that the fixed part of
\[\mathbf{R}^{\bullet}pr_{1,*}(\mathcal{H}om(\mathcal{I},
\mathcal{I})- \mathcal{O}_{C \times X})\] is just $T_C$. Since we
have already shown that the fixed loci are smooth, we can conclude
that $\Ob^f$ is zero. Therefore, by (\ref{eqn:DT}),
$N^{DT}_n(X,\beta)$ is equal to the sum of the integrals of
(\ref{equ:1}) and (\ref{equ:2}) over $C$. To do the integral, we
expand the fraction in terms of $[p\,]$ and integrate over $C$.
This proves the proposition. \qed

\begin{rem}
The method used here is a generalization of the method in \cite{MNOP1}. Similar ideas also appear in \cite{Song}.
\end{rem}

\bibliography{mainbiblio}
\bibliographystyle{plain}


\end{document}